\documentclass{amsart}
\usepackage{amsmath}
\usepackage{graphicx}
\usepackage{amsfonts}
\usepackage{amssymb}

\newtheorem{theorem}{Theorem}
\theoremstyle{plain}

\newtheorem{lemma}{Lemma}

\newtheorem{proposition}{Proposition}
\newtheorem{remark}{Remark}

\numberwithin{equation}{section}

\newlength{\originalbase}
\begin{document}
\title[The Mean Field Equation in a plane domain]
{The Mean Field Equation with Critical Parameter in a Plane Domain}
\author[Yilong Ni]{Yilong Ni\\
Department of Mathematics\\
University of Oklahoma\\
Norman, OK 73019}
\begin{abstract}
Consider the mean field equation with parameter $\lambda=8\pi$
in a bounded smooth domain $\Omega$. Denote by $E_{8\pi}(\Omega)$ the
infimum of the associated functional $I_{8\pi}(\Omega)$. We prove that if
$|\Omega|=\pi$, then $E_{8\pi}(\Omega)\ge E_{8\pi}(B_1)$ and equality
holds if and only if $\Omega$ is a ball. We also give a sufficient
condition for the existence of a minimizer for $I_{8\pi}(\Omega)$.
\end{abstract}
\maketitle
\section{Introduction}
Let $\Omega\subset{\mathbf R}^2$ be a bounded domain with
smooth boundary. Consider the following Mean Field equation
\begin{equation}
\label{eq}
\begin{cases}\displaystyle
&-\Delta  u=\frac{\lambda e^u}{\int_\Omega
e^u}, \mbox{ in }
\Omega,\\
&u|_{\partial\Omega}=0,
\end{cases}
\end{equation}
where
$\Delta=\frac{\partial^2}{\partial x_1^2}+\frac{\partial^2}{\partial x_2^2}$
and $\lambda$ is a real parameter. Equation (\ref{eq}) appears naturally in
many physical problems. For example in \cite{clmp} and \cite{clmp1},
it has been derived from the mean field limit of the Gibbs measure associated
to a system of $N$ vortices. It also arises in the study of the
Chern-Simons-Higgs model of superconductivity(see for example \cite{djlw1}).
To study the existence of solutions to eqaution (\ref{eq}), we may use the
variational approach. We consider the associated nonlinear functional
$I_\lambda$:
$$
I_\lambda(u,\Omega)=\frac{1}{2\lambda}\int_\Omega|\nabla
u|^2-\ln\left(\frac1{|\Omega|}\int_\Omega e^u\right),
$$
for $u\in H_0^1(\Omega)$, and denote
$$
E_\lambda(\Omega)=\inf_{u\in H^1_0(\Omega)}I_\lambda(u,\Omega).
$$
A well known fact is that
$I_\lambda(u)$ is bounded below if and only if $\lambda\le 8\pi$.
In particular, when $\lambda<8\pi$ Moser-Trudinger inequality \cite{mosert}
implies that the infimum of $I_\lambda(u,\Omega)$ is always
attained. However, in the critical case $\lambda=8\pi$,
the existence of a minimizer of $I_{8\pi}$ is a very difficult
problem and depends on the geometry of $\Omega$.
When $\Omega$ is a ball, the infimum of $I_{8\pi}$ is never
attained(see for example \cite{clmp}, \cite{cc}). Yet, when
$\Omega$ is thin, the infimum of $I_{8\pi}(\Omega)$ can be achieved(see for
example Proposition \ref{p1}). For general domains, there are only a few
results about the esitence of a minimizer of $I_{8\pi}(\Omega)$.
For example, Chang, Chen and Lin proved that the set of domains $\Omega$
on which the infimum of $I_{8\pi}$ is attained is open in the
$C^1$ topology(\cite{ccl}).
In this paper we study the functional $I_\lambda$ in the critical case
$\lambda=8\pi$ and obtain the following
\begin{theorem}\label{t1}
Suppose that $|\Omega|=|B_1|=\pi$, where $|\Omega|$ is the area
of $\Omega$ and $B_1$ is the unit ball in ${\mathbf R}^2$. Then
$$
E_{8\pi}(\Omega)\ge E_{8\pi}(B_1),
$$
and the equality holds if and only if $\Omega=B_1$.
\end{theorem}
If we view the infimum $E_{8\pi}(\Omega)$ as a Liouville type
energy of the domain, then Theorem \ref{t1} says that we can use
this energy to distinguish the unit ball from other domains with the
same area, since the unit ball has the lowest energy among these domains.
It is interesting to compare these with the Yamabe problem. Let $(M, g_0)$ be
a Riemannian manifold of dimension $n>2$. The Yamabe problem is to
find a metric $g$ conformal to $g_0$ such that $(M,g)$ has
constant scalar curvature $R$. If we write $g=u^{q-2}g_0$ with
$q=2n/(n-2)$, then $u$ satisfies the Yamabe equation:
$$
-L_{g_0} u= R u^{\frac{n+2}{n-2}},
$$
where $L_{g_0}$ is the conformal Laplacian. The associated
variational problem is
\begin{equation}
\label{mu} \mu(M,g_0)=\inf\left\{\int_M (a_n|du|^2+Ru^2)dV_{g_0} :
\int_M |u|^qdV_{g_0}=1 \right\},
\end{equation}
where $a_n=4(n-1)/(n-2)$. Denote the scalar curvature of $g$ by $R_{g}$,
we have
$$
\mu(M,g_0)=\inf_{g\in [g_0]}\frac{\int_M R_gd V_g}{\int_M dV_g},
$$
where $[g_0]$ is the conformal class of $g_0$. The solution to the famous
Yamabe problem can be summarized as the following theorem:
\begin{theorem}\label{t2}$($Yamabe, Trudinger and Aubin$)$
Let $(M,g_0)$ be a compact Riemannian manifold of dimension $n>2$.
Then $\mu(M)\le \mu(S^n)$, where $S^m$ is the sphere with the
standard metric. If $\mu(M)<\mu(S^n)$, the infimum of (\ref{mu})
is attained by a positive $C^\infty$ solution to the Yamabe
equation.
\end{theorem}
Based on the work of Yamabe and Trudinger, Aubin and
Schoen \cite{schoen} completed the solution of the Yamabe problem
by proving that $\mu(M,g_0)< \mu(S^n)$ unless $M$ is conformally equivalent
to the standard sphere. This indicates that we can use $\mu(M,g_0)$
to distinguish the standard sphere from other compact manifolds in conformal
sense. We are certainly wondering whether $E (\Omega)$ plays a similar
role in the mean field equation as $\mu(M,g_0)$ in the Yamabe problem.
Note that if $\mu(M)<\mu(S^n)$, then
the infimum of (\ref{mu}) is attained by a positive $C^\infty$
solution to the Yamabe equation. It is interesting to consider,
for a bounded smooth domain $\Omega\subset{\mathbf R}^2$ satisfying
$|\Omega|=\pi$ and $E_{8\pi}(\Omega)>E_{8\pi}(B_1)$, when
the infimum of $I(\cdot,\Omega)$ is attained by a smooth solution to the
mean field equation. Chang, Chen and Lin \cite{ccl} gave an example of a
dumbell $\Omega_h$ which consists of two disjount balls $B(r_1)$ and
$B(r_2)$ connected with a tube of small width $h>0$. They proved that when
$r_1<r_2$ and $h$ is sufficiently small, the infimum of $I(\cdot,\Omega_h)$
is not attained(see Proposition 7.3 in \cite{ccl}). Therefore in order that
the infimum of $I(\cdot,\Omega)$ is attained, we must add more conditions
on the domain $\Omega$. For example, we could require that $\Omega$ is thin.
\begin{proposition}\label{p1}
Suppose $|\Omega|=\pi$ and $\Omega$ can be covered by a strip with
width $d\le \frac{\pi}{2\sqrt{e}}=0.9527...$, then the infimum of
$I_{8\pi}(u,\Omega)$ can be achieved by a function in
$H_0^1(\Omega)$.
\end{proposition}
The paper is organized as follows. In section 2, we study the
regular part $\gamma(\Omega)$ of the Green's function of $\Omega$.
The property of $\gamma(\Omega)$ will be used in the proof of
Theorem \ref{t1}. In section 3 we derive some standard estimates
for $E_{8\pi}(\Omega)$. Theorem \ref{t1} then follows from these estimates.
The existence result Proposition \ref{p1} is proved in Section 4.

We will suppress the subscript ``$8\pi$" if no confusion would result.
\medskip
\section{Regular Part of the Green's Function}
Denote by $G(x,y)$ the Green's function of $\Omega$:
$$
\begin{cases}\displaystyle &-\Delta_x G(x,y)=\delta_y(x), \mbox{ in } \Omega,\\
&G(\cdot,y)|_{\partial\Omega}=0.
\end{cases}
$$
Let
$$
\gamma(x,y)=G(x,y)-\frac1{2\pi}\ln\frac1{|x-y|},
$$
be the regular part of $G$, and set $\gamma(x)=\gamma(x,x)$,
$\gamma(\Omega)=\sup_{x\in\Omega}(\gamma(x))$. Then we have the
following Lemma.
\begin{lemma}
\label{gamma0} Suppose $\Omega\subset{\mathbf R}^2$ is an bounded
domain with $|\Omega|=\pi$. Then $\gamma(\Omega)\le 0$, and
equality holds if and only if $\Omega=B_1$, the unit ball in
${\mathbf R}^2$.
\end{lemma}
\begin{proof}
For any $x_0\in\Omega$,
$$
G(x_0,y)=\frac1{2\pi}\log\frac1{|x_0-y|}+\gamma(x_0)+O(|x_0-y|),
\mbox{ as }y\to x_0.
$$
Therefore for any $\epsilon>0$, there exists $\rho>0$, s.t.
$$
\left|G(x_0,y)-\frac1{2\pi}\log\frac1{|x_0-y|}-\gamma(x_0)\right|
\le\epsilon\mbox{ whenever }|y-x_0|\le\rho.
$$
It follows that when $\tau$ is sufficiently large,
\begin{equation}
\label{31}
B_{x_0}(e^{-2\pi(\tau-\gamma(x_0)+\epsilon)})
\subset\Omega_\tau\subset
 B_{x_0}(e^{-2\pi(\tau-\gamma(x_0)-\epsilon)}),
\end{equation}
where $\Omega_\tau= \{G(x_0,y)>\tau\}$ and $B_{x_0}(r)$ is the ball in
${\mathbf R}^2$ with radius $r$ and centered at $x_0$.

Let $G^*(y):\ B(1)\to \mathbf R$ be the rearrangement of $G(x_0,y)$
and let
$$
\mu(\tau)=|\Omega_\tau|,\quad \rho(\tau)=\sqrt{\mu(\tau)/\pi}.
$$
Then (\ref{31}) implies
\begin{equation}
\label{32} e^{-2\pi(\tau-\gamma(x_0)+\epsilon)}\le \rho(\tau) \le
e^{-2\pi(\tau-\gamma(x_0)-\epsilon)}.
\end{equation}
Define a function $\phi(t): [0,\infty)\to\mathbf R$ by
$$
\phi(t)=\int_{\Omega\setminus\Omega_t}|\nabla_y G(x_0,y)|^2dV_y
-\int_{B(1)\setminus B(\rho(t))}|\nabla G^*(y)|^2dV_y.
$$
From the properties of rearrangement, it is easy to see that $\phi(t)$
is increasing, $\phi(t)\ge0$, $G^*(y)|_{\partial B(\rho(t))}=t$ and
$G^*(y)|_{\partial B(1)}=0$. It follows that
$$
\int_{B(1)\setminus B(\rho(t))}|\nabla G^*(y)|^2dV_y\ge
\int_{B(1)\setminus B(\rho(t))}|\nabla G_0(y)|^2dV_y,
$$
where $G_0(y)=t\log(|y|)/\log(\rho(t))$. Therefore
using (\ref{32}), for sufficiently large $t$,
\begin{align}
\label{33}
\phi(t)=&\int_{\Omega\setminus \Omega_t}
|\nabla_y G(x_0,y)|^2dV_y
-\int_{B(1)\setminus B(\rho(t))}|\nabla G^*(y)|^2dV_y\\
\nonumber
\le & \int_{\Omega\setminus\Omega_t}\nabla_y G(x_0,y)
\cdot\nabla_y G(x_0,y)dV_y
- \int_{B(1)\setminus B(\rho(t))}|\nabla G_0(y)|^2dV_y\\
\nonumber
=&-t\int_{\partial\Omega_t} \frac{\partial
G(x_0,y)}{\partial{\mathbf n}}dS
 + \frac{2\pi t^2} {\log\rho(t)}\\
\nonumber
\le & t-\frac{t^2}{t-\gamma(x_0)+\epsilon}
=(\epsilon-\gamma(x_0))\frac{t}{t-\gamma(x_0)+\epsilon}.
\end{align}
Since $\phi(t)\ge0$ and $\epsilon$ is arbitraty, we obtain that
$\gamma(x_0)\le 0$. Hence $\gamma(\Omega)\le 0$.
Furthermore, if $\gamma(\Omega)= 0$, it follows from the continuity
of $\gamma(x)$ that there exists $x_0$, s.t. $\gamma(x_0)=0$. Then
(\ref{33}) implies $\phi(t)\le\epsilon$. Since $\phi(t)$ is increasing,
$\phi(t)\ge0$ and $\epsilon$ is arbitrary, we have $\phi(t)=0$, for
all $t>0$. From the properties of rearrangement, we have $\Omega=B(1)$.
\end{proof}
\begin{remark}
The proof of Lemma \ref{gamma0} can be easily extended to higher
dimensions. Let $\Omega\subset{\mathbf R}^n$, $n\ge 2$ and let
$G_n(x,y)$ be the the solution of the following differential
equation:
$$
\begin{cases}\displaystyle &-\Delta_n u(x)=\delta_y(x), \mbox{ in } \Omega,\\
&u|_{\partial\Omega}=0,
\end{cases}
$$
where $\Delta_n u(x)=\nabla(|\nabla u(x)|^{n-2}\nabla u(x))$ is the
$n-$Laplacian. We may consider as in the case $n=2$ the regular part:
$$
\gamma_n(x,y)=G_n(x,y)-\omega_{n-1}^{-\frac1{(n-1)}}\ln\frac1{|x-y|},
$$
where $\omega_{n-1}$ is the volume of the $(n-1)$-dimensional
sphere. Set $\gamma_n(x)=\gamma_n(x,x)$, and
$\gamma(\Omega)=\sup_{x\in\Omega}(\gamma_(x))$. Then we have the
following Lemma.
\begin{lemma}
\label{gamman} Suppose $\Omega\subset{\mathbf R}^n$ is an bounded
domain with $|\Omega|=|B_1^n|$, where $B_1^n$ the unit ball in
${\mathbf R}^n$. Then $\gamma_n(\Omega)\le 0$, and equality holds
if and only if $\Omega=B_1^n$.
\end{lemma}
\end{remark}
\medskip
\section{Estimates of $E(\Omega)$}
In this section, we always assume that $\Omega\subset{\mathbf R}^2$
is a bounded domain with smooth boundary and $|\Omega|=\pi$.
We have the following estimates for $E_{8\pi}(\Omega)$:
\begin{proposition}\label{l1}
$$
E_{8\pi}(\Omega)\le -1 -4\pi\gamma(\Omega).
$$
\end{proposition}
\begin{proposition}\label{l2}
If the infimum of $I_{8\pi} (\cdot,\Omega)$ is not attained, then
$$
E_{8\pi}(\Omega)\ge -1-4\pi\gamma(\Omega).
$$
\end{proposition}
These two estimates are very standard. For example, Propostion \ref{l2} has
been proved in \cite{clmp}. For the sake of completeness, we shall reprove
these two estimates using a different method.

Suppose $G(x,y)$ is the Green's function of $\Omega$, and
$\gamma(x)$ is its regular part. We may assume that $\gamma(x)$
attains its maximum value at $x_0\in \Omega$. Without loss of
generality, suppose $x_0=0$. Let
\begin{equation}
\label{gg} G(x)=8\pi G(x,0) =4\ln\frac1{|x|}+A+\alpha(x),
\end{equation}
where $A=8\pi\gamma(0)$. For any $\Lambda,\epsilon>0$, such that
$\rho=\Lambda\epsilon<\mbox{dist}(0,\partial\Omega)$, we choose a
test function
$$
\phi(x)=\begin{cases}\displaystyle
&2\ln\frac{1}{\epsilon^{2}+|x|^{2}} -
C_\epsilon,\hskip 1cm |x|\le \rho\\
&G(x)-\eta(x)\alpha(x),\hskip 1cm  \rho\le|x|\le 2\rho\\
&G(x),\hskip 3cm |x|\ge 2\rho
\end{cases}
$$
where $\eta(x)$ is a $C^\infty$ bump function:
$$
\eta(x)=\begin{cases}\displaystyle &1, \ |x|\le \rho\\
&0, \ |x|\ge 2\rho
\end{cases}
$$
satisfying $|\nabla\eta(x)|\le 2/\rho$, and
$$
C_\epsilon=2\ln\frac{\Lambda^{2}}{1+\Lambda^{2}}-A
$$
so that the function $\phi(x)\in H_0^{1}(\Omega)$. It suffices to
prove that
$$
I(\phi,\Omega)\le -1-\frac{A}{2}.
$$
First
$$
\int_{B(\rho)} |\nabla \phi|^2= 4\int_{B(\rho)}
\left(\frac{2r}{\epsilon^2+r^2}\right)^2dx\\
=  32\pi\int_0^\rho\frac{r^2r}
{(\epsilon^{2}+r^{2})^2}dr\\
=  16\pi\int_{ \sigma }^1\frac{(1-t)^2}{t}dt,
$$
where
$$
t=\frac{\epsilon^{2}}{\epsilon^{2}+r^{2}}\mbox{ and }
\sigma=t(\rho)=\frac1{1+\Lambda^2}.
$$
Therefore
\begin{equation}
\label{brho} \frac1{16\pi}\int_{B(\rho)} |\nabla
\phi|^2=-\ln\sigma-1+\sigma.
\end{equation}
Second, since $G(x)=-4\ln |x|+A+\alpha(x)$, $\alpha(x)$ is a
smooth function and $\alpha(0)=0$, we have when $\rho$ is
sufficiently small,
\begin{align*}
\int_{\Omega\setminus B(\rho)}|\nabla G|^2= & -\int_{\partial
B(\rho)}G\frac{\partial G}{\partial {\mathbf n}}= -\left(
4\ln\frac1\rho+A+O(\rho) \right)\int_{\partial
B(\rho)}\frac{\partial G}{\partial {\mathbf n}} \\
=&32\pi \ln\frac1\rho + 8\pi A+O(\rho).
\end{align*}
For $|x|>\rho$, we have $\phi(x)=G(x)-\eta(x)\alpha(x)$. Therefore
$$
|\nabla \phi|^2=|\nabla G-\nabla(\eta\alpha)|^2 = |\nabla
G|^2+2\nabla G\cdot\nabla (\eta\alpha)+|\nabla (\eta\alpha)|^2.
$$
Since $\alpha(x)$ is smooth and $\alpha(0)=0$, $|\nabla
G|=\frac2\rho(1+O(\rho))$ for $2\rho>|x|>\rho$, and
$|\nabla(\eta\alpha)|\le C$ for any $x$. Hence,
$$
\int_{\Omega\setminus B(\rho)}2\nabla G\cdot\nabla
(\eta\alpha)+|\nabla (\eta\alpha)|^2 =\int_{B(2\rho)\setminus
B(\rho)}2\nabla G\cdot\nabla (\eta\alpha)+|\nabla
(\eta\alpha)|^2=O(\rho),
$$
which implies
\begin{equation}
\label{brhoc} \int_{\Omega\setminus B(\rho)}
|\nabla\phi|^2=32\pi\ln\frac1\rho + 8\pi A+O(\rho).
\end{equation}
Combining (\ref{brho}) and (\ref{brhoc}), we obtain
\begin{equation}
\label{item1} \frac1{16\pi}\int_{\Omega} |\nabla
\phi|^2=-\ln\sigma-1+\sigma -2\ln\rho+\frac A 2+O(\rho).
\end{equation}

Now we turn to the estimate of $\int_\Omega e^\phi$.
\begin{align*}
\int_{B(\rho)} e^\phi=&
\exp(-C_\epsilon)\int_{B(\rho)} \left(
\frac{1}{\epsilon^{2}+r^{2}} \right)^2dx\\
=&2\pi e^{-C_\epsilon}\int_0^\rho
\frac{r}{\epsilon^{2}+r^{2}}dr
=\frac{\pi e^A}{\epsilon^2}\frac{1+\Lambda^2}{\Lambda^2}
\end{align*}
Since $\phi(x)>0$ when $\rho$ is small, we have
\begin{equation}\label{item2}
-\ln\frac{1}{\pi}\int_\Omega e^\phi\le
-\ln\frac{1}{\pi}\int_{B(\rho)} e^\phi=-A+2\ln\epsilon- \ln\frac{1
+\Lambda^2}{\Lambda^2}.
\end{equation}
Sending $\Lambda\to\infty$ and using (\ref{item1}) and
(\ref{item2}) we finally get
\begin{align*}
I(\phi,\Omega)=&\frac1{16\pi}\int_\Omega|\nabla
\phi|^2-\ln\frac{1}{\pi}\int_\Omega e^\phi\\
\le & -1-\frac A2=-1-4\pi\gamma(\Omega)
\end{align*}
This completes the proof of Proposition \ref{l1}.

\medskip
Next, we prove the opposite inequality, Proposition \ref{l2}. We will use a
similar argument given in \cite{djlw}. Suppose the infimum of
$I (\cdot,\Omega)$ is not attained. For $\epsilon>0$ we define
$$
I^\epsilon(u)=\frac{1}{16\pi}\int_\Omega|\nabla
u|^2-(1-\epsilon)\ln\frac{1}{\pi}\int_\Omega e^u,
$$
and
$$
E^\epsilon=\inf_{u\in H^{1}_0(\Omega)}I^\epsilon(u).
$$
\begin{theorem}
$($\cite{zhu}$)$ Let $\Omega\subset {\mathbf R}^2$ be a bounded smooth domain,
$\Omega^*$ be the ball in ${\mathbf R}^2$ which has the same area as $\Omega$,
and denote
$$
D_{a,b}(\Omega)=\{f(x)-b\in H^1_0(\Omega)\ : \ \int_\Omega e^{2f}dx=a \}.
$$
We have the following sharp inequality:
$$
\inf_{w\in D_{a,b}(\Omega)}\int_\Omega|\nabla w|^2dx\ge 4\pi
(\ln\frac{ae^{-2b}}{\pi r^2}+\frac{\pi r^2}{ae^{-2b}}-1),
$$
where $r$ is the radius of $\Omega^*$.
\end{theorem}
It follows from the above sharp inequality that for any $u\in H^1_0(\Omega)$
\begin{equation}\label{bddb}
\int_\Omega e^udx\le \pi e\exp\left(\frac{1}{16\pi}
\int_\Omega |\nabla u|^2dx\right),
\end{equation}
which implies $I^\epsilon(u)$ is bounded below by $-1$.
\begin{lemma}
$E^\epsilon$ is achieved by a function $u_\epsilon \in H^{1}_0(\Omega)$,
which is the solution of the following equation
\begin{equation}\label{eqne}
 \begin{cases}\displaystyle &-\Delta u_\epsilon=\frac{(1-\epsilon)8\pi e^{u_\epsilon}}{\int_\Omega
e^{u_\epsilon}}, \mbox{ in }
\Omega,\\
&u_\epsilon|_{\partial\Omega}=0,
\end{cases}
\end{equation}
\end{lemma}
\begin{proof}
For fixed $\epsilon>0$, let $\{u_n\}$ be a minimizing sequence of
$\inf I^\epsilon$. If follows from (\ref{bddb}) that
\begin{align*}
I^\epsilon(u_n)=&\frac{1}{16\pi}\int_\Omega|\nabla
u_n|^2-(1-\epsilon)\ln\frac{1}{\pi}\int_\Omega e^{u_n}\\
=&(1-\epsilon)I(u_n)+\epsilon\int_\Omega|\nabla u_n|^2\\
\ge & -(1-\epsilon)+\epsilon\int_\Omega|\nabla u_n|^2.
\end{align*}
Therefore $||\nabla u_n||_{L^2}<C$. Since $u_n|_{\partial\Omega}=0$,
using Poincar\'e's inequality we obtain that $||u_n||_{H^1}<C$. Hence
$u_n\rightharpoonup u_\epsilon$ in $H^1_0(\Omega)$ for some
$u_\epsilon\in H^1_0(\Omega)$. It follows from Trudinger's inequality
\cite{trudinger} that $e^{ku_n}\to e^{ku_\epsilon}$ in $L^1(\Omega)$ for any
$k>0$. Therefore $u_\epsilon$ is a minimizer and satisfies the Euler-Lagrange
equation (\ref{eqne}).
\end{proof}
It follows from standard elliptic estimates that $u_\epsilon\in
C^\infty(\Omega)$. Suppose $u_\epsilon$ attains its maximum at
$x_\epsilon\in\Omega$ and set $\lambda_\epsilon=u_\epsilon(x_\epsilon)
=\max_{x\in\bar{\Omega}}u_\epsilon(x)$. The following Lemma is immediate.
\begin{lemma}
\label{l21}
There exists a subsequence $\epsilon_i\to 0$, such that
$$
\lim_{i\to \infty}\lambda_{\epsilon_i}=+\infty.
$$
\end{lemma}
\begin{proof}
Suppose $\lambda_\epsilon$ is bounded above, $\lambda_\epsilon\le
C<+\infty$, then
$$
\int_\Omega e^{u_\epsilon}\le C.
$$
Since $E^\epsilon\le I^\epsilon(0)=0$, we have $\int_\Omega|\nabla
u_\epsilon|^{2}\le C$. Hence there exists a
subsequence of $u_\epsilon$ which converge weakly to $u_0$ in
$H^{1}_0(\Omega)$. We can easily check that $u_0$ is a minimizer
for $I(\cdot,\Omega)$, which contradict with the assumption that
the infimum of $I(\cdot,\Omega)$ is not attained.
\end{proof}
In the following, for simplicity, we shall not distinguish a subsequence
$\{\epsilon_i\}$ from the original $\{\epsilon\}$.

Next, we claim that $x_\epsilon$ will stay away from $\partial\Omega$, which
implies
$$
x_\epsilon\to \bar{x}\in\Omega\mbox{ as }\lambda\to 0.
$$
The claim can be proved
by the moving plane method(see \cite{gnn}) and an interior integral
estimate(cf. page 163-164 in \cite{han}). We shall omit the details here.
Set $\tau_\epsilon=e^{\lambda_\epsilon/2}$ and
$$
\alpha_\epsilon=\left(\frac{(1-\epsilon)\pi}{
\int_\Omega e^{u_\epsilon}}\right)^{1/2}\tau_\epsilon
$$
If $\alpha_\epsilon$ stays bounded as $\epsilon\to 0$, standard elliptic
estimate of (\ref{eqne}) implies that $u_\epsilon$ is uniformly bounded
as $\epsilon\to 0$, which contradicts with the fact that $u_\epsilon$ blows
up(Lemma \ref{l21}). Therefore we have:
\begin{lemma}
\label{l22}
$$
\lim_{\epsilon\to 0}\alpha_\epsilon=+\infty.
$$
\end{lemma}
\noindent Define
\begin{equation}\label{phie}
\phi_\epsilon(x)=u_\epsilon(\alpha_\epsilon^{-1}x+x_\epsilon)-2\ln\tau_\epsilon.
\end{equation}
We can easily see that $\phi_\epsilon$ satisfies
\begin{equation}
\begin{cases}\displaystyle &-\Delta \phi_\epsilon=8e^{\phi_\epsilon},\mbox{ in
}\Omega_\epsilon\\
&\phi_\epsilon|_{\partial\Omega_\epsilon}=-2\ln\tau_\epsilon,
\end{cases}
\end{equation}
where $\Omega_\epsilon=\alpha_\epsilon\cdot(\Omega-x_\epsilon)$.
We claim that for any $R>0$, $\phi_\epsilon$ is bounded in $B(R)$
uniformly in $\epsilon$. In fact, let $\phi_\epsilon^{(1)}$ be the
unique solution to
$$
\begin{cases}\displaystyle &-\Delta \phi_\epsilon^{(1)}=8 e^{\phi_\epsilon}, \mbox{
in }B(2R)\\
&\phi^{(1)}_\epsilon|_{\partial B(2R)}=0.
\end{cases}
$$
Since $x_\epsilon$ is a maximum point of $u_\epsilon(x)$, we
have $\phi_\epsilon\le \phi_\epsilon(0)=0$ and $e^{\phi_\epsilon}\le 1$. It
follows that $||\phi_\epsilon^{(1)}||_{L^\infty}\le C<+\infty$.
Let $\phi_\epsilon^{(2)}=\phi_\epsilon-\phi_\epsilon^{(1)}$. Then
$\phi_\epsilon^{(2)}\le -\phi_\epsilon^{(1)}\le C$. Since
$2C-\phi_\epsilon^{(2)}\ge C$ is harmonic, positive and
$$
2C-\phi_\epsilon^{(2)}(0)=2C-\phi_\epsilon(0)+\phi^{(1)}_\epsilon(0)
 \le 3C,
$$
Harnack's inequality implies that
$||2C-\phi_\epsilon^{(2)}||_{L^\infty}\le \tilde{C}$ in $B(R)$. Hence
$||\phi_\epsilon ||_{L^\infty(R)}\le C+\tilde{C}$. Elliptic
estimates yield that, up to a subsequence, $\phi_\epsilon \to
\phi_0$ in $C^{1,\alpha}(B(R/2))$ for some $\alpha\in(0,1)$ and
$\phi_0$ satisfies
$$
\begin{cases}\displaystyle &-\Delta \phi_0=8 e^{\phi_0}, \mbox{
in }{\mathbf R}^2\\
&\phi_0(0)=0. \end{cases}
$$
Since $\int_{\Omega_\epsilon}e^{\phi_\epsilon}=8\pi(1-\epsilon)$,
we get
$$
\int_{{\mathbf R}^2}e^{\phi_0}\le \overline{\lim_{R\to\infty}}\,
\overline{\lim_{n\to\infty}} \int_{B(R)}e^{\phi_{\epsilon_n}}\le
8\pi(1-\epsilon).
$$
The uniqueness theorem in \cite{chenli} implies that
\begin{equation}
\label{phi0}
\phi_0(x)=2\ln\frac{1}{1+|x|^2}.
\end{equation}
The following Lemma is due to Brezis and Merle\cite{bm}:
\begin{lemma}\label{bm}
Let $D\subset {\mathbf R}^2$ be a bounded domain and $u$ be a solution to the
following equation,
$$
 \begin{cases}\displaystyle &-\Delta u=f(x), \mbox{ in } D,\\
&u|_{\partial D}=0.
\end{cases}
$$
If $f\in L^1(D)$, then for any $\delta\in(0,4\pi)$, there is a constant
$C(\delta)$ such that
$$
\int_D \exp\left(\frac{(4\pi-\delta)|u(x)|}{||f||_{L^1(D)}}\right)\le
C(\delta).
$$
\end{lemma}
\noindent Using this lemma we have:
\begin{lemma}
\label{l23}
For any $K\subset\subset \Omega\setminus\{\bar{x}\}$, there exists a constant
$C(K)$, such that $u_\epsilon(x)\le C(K)$, for all $x\in K$.
\end{lemma}
\begin{proof}
From (\ref{eqne}) and Lemma \ref{bm}, we know that $e^{u_\epsilon}
\in L^p(\Omega)$ for $p\in(0,\frac12)$. For any given
$K\subset\subset\Omega\setminus\{\bar{x}\}$, we choose $K'$ such that
$K\subset\subset K'\subset\subset\Omega\setminus\{\bar{x}\}$.
Since $\phi_\epsilon\to\phi_0$ in $C^{1,\alpha}$ and
$\int_{{\mathbf R}^2}e^{\phi_0}=\pi$, we obtain that
\begin{equation}\label{to0}
\lim_{\epsilon\to 0}\frac{\int_{K'}e^{u_\epsilon}}
{\int_{\Omega} e^{u_\epsilon}}=0.
\end{equation}
Let $u_\epsilon^{(1)}$ be the unique solution to
\begin{equation}\label{u1}
 \begin{cases}\displaystyle &-\Delta u_\epsilon^{(1)}=\frac{(1-\epsilon)8\pi e^{u_\epsilon}}
{\int_\Omega e^{u_\epsilon}}, \mbox{ in } K', \\
&u_\epsilon^{(1)}|_{\partial\Omega}=0,
\end{cases}
\end{equation}
It follows from (\ref{to0}) and Lemma \ref{bm} that $e^{u_\epsilon^{(1)}}
\in L^p(K')$ for some $p>1$. Since $u_\epsilon^{(2)}:= u_\epsilon
-u_\epsilon^{(1)}$ is harmonic in $K'$, Harnack's inequality implies that
\begin{align*}
||u_\epsilon^{(2)}||_{L^\infty(K')}&\le C||u_\epsilon^{(2)}||_{L^p(K')}\\
&\le C(||u_\epsilon||_{L^p(K')}+||u_\epsilon^{(1)}||_{L^p(K')})\\
&\le C.
\end{align*}
Therefore, for some $p>1$ we have
$$
\int_{K'}e^{pu_\epsilon}=\int_{K'}e^{pu_\epsilon^{(1)}}
\cdot e^{pu_\epsilon^{(2)}}\le C\int_{K'}e^{pu_\epsilon^{(1)}}\le C
$$
It follows from the standard elliptic estimates of (\ref{u1}) that
$||u_\epsilon^{(1)}||_{L^\infty(K)}\le C$,
thus $||u_\epsilon||_{L^\infty(K)}\le C$.
\end{proof}
We may assume without loss of generality that $\bar{x}=0$.
Since $u_\epsilon$ satisfies (\ref{eqne}), it follows from Lemma
\ref{l23} that $u_\epsilon\to G(x)$ in
$C^{1,\alpha}_{loc}(\Omega\setminus \{\bar{x}\})$, where $G(x)$ was defined in
(\ref{gg}).
\begin{lemma}
\label{l24} For fixed $R$, let $r_\epsilon=R/\alpha_\epsilon$.
Then for any $x\in \Omega\setminus B(r_\epsilon)$,
$$
u_\epsilon(x)\ge G(x)-\lambda_\epsilon+2\ln\left(
\frac{1}{\pi}\int_\Omega e^{u_\epsilon}\right)+2\ln\frac{R^{2}}
{1+R^{2}}-A +o_\epsilon(1),
$$
where $o_\epsilon(1)$ stands for some function that goes to $0$ as
$\epsilon\to 0$.
\end{lemma}
\begin{proof}
On $\partial B(r_\epsilon)$, $G(x)$ and $u_\epsilon$ have the following
asymptotic behavior:
\begin{align*}
G(x)&=-4\ln r_\epsilon+A+o_\epsilon(1),\\
u_\epsilon(x)&=\lambda_\epsilon+2\ln\frac{1}{1+R^{2}}+o_\epsilon(1).
\end{align*}
Therefore on $\partial B(r_\epsilon)$,
\begin{align*}
u_\epsilon-G  =&\lambda_\epsilon+2\ln\frac{1}{1+R^{2}}+
4\ln \frac{R}{\alpha_\epsilon}- A +o_\epsilon(1)\\
=&\lambda_\epsilon+2\ln\frac{R^{2}}{1+R^{2}}-A-
\ln\frac{(1-\epsilon)\pi \tau_\epsilon^2
}{\int_\Omega e^{u_\epsilon}} +o_\epsilon(1)\\
=& -\lambda_\epsilon +2\ln\left( \frac{1}{\pi}\int_\Omega
e^{u_\epsilon}\right)+2\ln\frac{R^{2}} {1+R^{2}}-A +o_\epsilon(1).
\end{align*}
Let
\begin{equation}\label{depsilon}
D_\epsilon =- \lambda_\epsilon +2\ln\left( \frac{1}{\pi}\int_\Omega
e^{u_\epsilon}\right)+2\ln\frac{R^{2}} {1+R^{2}}-A,
\end{equation}
and consider the function
$u_\epsilon -G- D_\epsilon$ on $\Omega\setminus B(r_\epsilon)$. It
satisfies
$$
\begin{cases}\displaystyle &\Delta (u_\epsilon -G- D_\epsilon)\le 0,\mbox{ in }
\Omega\setminus B(r_\epsilon)\\
&u_\epsilon -G- D_\epsilon\ge o_\epsilon(1),\mbox{ on }\partial
(\Omega\setminus B(r_\epsilon)). \end{cases}
$$
Then the lemma follows immediately form the maximum principle.
\end{proof}

Now we estimate $I(u_\epsilon)$. For fixed $R$, let $r_\epsilon
=R/\alpha_\epsilon$ and choose $\delta>r_\epsilon$, such
that $B(\delta)\subset \Omega$. Then
$$
\int_\Omega|\nabla u_\epsilon|^2=\left(\int_{\Omega\setminus
B(\delta)} +\int_{B(\delta)\setminus B(\epsilon)}
+\int_{B(\epsilon)}  \right)|\nabla u_\epsilon|^2:= {\mathbf
I}_1+{\mathbf I}_2+{\mathbf I}_3.
$$
Since $\phi_\epsilon\to\phi_0$ in $C^{1,\alpha}(B(R/2))$
(\ref{phi0}), we obtain that
\begin{align*}
{\mathbf I}_3=&\int_{B(r_\epsilon)}|\nabla u_\epsilon|^2=
\int_{B(R)}\left|\nabla  \left(2\ln\frac{1}{1+|x|^{2}}\right)
\right|^2+o_\epsilon(1)\\
=&32\pi\int_0^R\frac{r^{2}r}{(1+r^{2})^2}dr+o_\epsilon(1)
=16\pi\left( \ln(1+R^{2})- \frac{R^2}{1+R^2} \right)
+o_\epsilon(1).
\end{align*}
From $u_\epsilon\to G(x)$ in $C^{2}_{loc}(\Omega\setminus\{0\})$, we know that
$$
{\mathbf I}_1= \int_{\Omega\setminus B(\delta)}|\nabla
u_\epsilon|^2=\int_{\Omega\setminus B(\delta)}|\nabla
G|^2+o_\epsilon(1) =  -\int_{\partial  B(\delta)} G\frac{\partial
G}{\partial{\mathbf n}}+o_\epsilon(1).
$$
To estimate ${\mathbf I}_2$, we apply Lemma \ref{l24}
\begin{align*}
{\mathbf I}_2=&-\int_{ B(\delta)\setminus B(r_\epsilon)}\Delta
u_\epsilon\cdot u_\epsilon+ \int_{\partial(B(\delta)\setminus
B(r_\epsilon))}  u_\epsilon
\frac{\partial u_\epsilon}{\partial{\mathbf n}} \\
\ge &\int_{ B(\delta)\setminus B(r_\epsilon)}(-\Delta
u_\epsilon)\cdot (G+D_\epsilon+o_\epsilon(1))+
\int_{\partial(B(\delta)\setminus B(r_\epsilon))}  u_\epsilon
\frac{\partial u_\epsilon}{\partial{\mathbf n}}\\
=&\int_{\partial(B(\delta)\setminus B(r_\epsilon))}\left(-\frac{\partial
u_\epsilon}{\partial{\mathbf n}}G+\frac{\partial G}{\partial{\mathbf n}}
u_\epsilon\right)+\int_{\partial(B(\delta)\setminus B(r_\epsilon))}(u_\epsilon
-D_\epsilon) \frac{\partial u_\epsilon}{\partial{\mathbf n}}+o_\epsilon(1).
\end{align*}
It follows that
\begin{align*}
\int_{\Omega\setminus B(r_\epsilon)}|\nabla u_\epsilon|^2\ge &
\int_{\partial B(\delta)}
\left(u_\epsilon \frac{\partial u_\epsilon}{\partial{\mathbf
n}} -G \frac{\partial G}{\partial{\mathbf n}}
-D_\epsilon\frac{\partial u_\epsilon}{\partial{\mathbf n}}
+u_\epsilon \frac{\partial G}{\partial{\mathbf
n}} -G \frac{\partial u_\epsilon}{\partial{\mathbf n}} \right)\\
&\quad-\int_{\partial B(r_\epsilon)}\left( \frac{\partial
G}{\partial{\mathbf n}}u_\epsilon + \frac{\partial
u_\epsilon}{\partial{\mathbf n}}
(u_\epsilon-G-D_\epsilon)\right)+o_\epsilon(1).
\end{align*}
Since $u_\epsilon\to G(x)$ in $C^{1,\alpha}_{loc}(\Omega\setminus
\{0\})$, we can easily see that
$$
\int_{\partial B(\delta)}u_\epsilon \frac{\partial
u_\epsilon}{\partial{\mathbf n}} -G \frac{\partial
G}{\partial{\mathbf n}}=o_\epsilon(1),\quad
\int_{\partial B(\delta)}u_\epsilon \frac{\partial G}
{\partial{\mathbf n}} -G \frac{\partial
u_\epsilon}{\partial{\mathbf n}}=o_\epsilon(1),
$$
and on $\partial B(r_\epsilon)$,
\begin{align*}
G(x)& =4\ln\frac1{|x|}+A+o_\epsilon(1),\\
\frac{\partial G}{\partial{\mathbf n}}&=-4\frac{1}
{r_\epsilon}+o_\epsilon(1),\\
u_\epsilon(x)&=\lambda_\epsilon+2\ln\frac{1}{1+R^{2}}+o_\epsilon(1),\\
\frac{\partial u_\epsilon}{\partial{\mathbf n}}&=-\left(\frac{4R
}{1+R^{2}}+o_\epsilon(1)\right)\alpha_\epsilon.
\end{align*}
It follows that
$$
\int_{B(r_\epsilon)}\frac{\partial G}{\partial{\mathbf
n}}u_\epsilon=(\lambda_\epsilon+2\ln\frac{1}{1+R^2}+o_\epsilon(1))
\int_{B(r_\epsilon)}\frac{\partial G}{\partial{\mathbf n}}= -8\pi
(\lambda_\epsilon+2\ln\frac{1}{1+R^2}+o_\epsilon(1)),
$$
and
$$
\int_{B(r_\epsilon)}\frac{\partial u_\epsilon}{\partial{\mathbf
n}}(u_\epsilon-G-D_\epsilon)=o_\epsilon(1).
$$
Since $u_\epsilon$ satisfies (\ref{eqne}), we have
$$
\int_{\partial (B(\delta))} \frac{\partial
u_\epsilon}{\partial{\mathbf n}}=\int_{B(\delta)} \Delta
u_\epsilon=\int_{B(\delta)}\frac{-(1-\epsilon)8\pi e^{u_\epsilon}}
{\int_\Omega e^{u_\epsilon}}\ge -(1-\epsilon)8\pi.
$$
Combining the estimates for $I_1$, $I_2$ and $I_3$, we finally have
\begin{align*}
\frac1{16\pi}\int_\Omega|\nabla u_\epsilon|^2\ge &
\ln(1+R^2)-\frac{R^2}{1+R^2}+\frac {D_\epsilon}2(1-\epsilon)
+\frac{\lambda_\epsilon}2+\ln\frac{1}{1+R^2} +o_\epsilon(1)\\
\ge &(1-\epsilon)\left(\ln\left(\frac1\pi\int_\Omega
e^{u_\epsilon}\right)+\ln\frac{R^2}{1+R^2}-\frac A2
\right)-\frac{R^2}{1+R^2}+o_\epsilon(1).
\end{align*}
Hence
\begin{align*}
E^\epsilon=&I^\epsilon(u_\epsilon)=\frac{1}{16\pi}\int_\Omega |\nabla
u_\epsilon|^2-(1-\epsilon)\ln\left(\frac1\pi\int_\Omega
e^{u_\epsilon}\right)\\
\ge& -\frac{(1-\epsilon)}2A+(1-\epsilon)\ln\frac{R^2}{1+R^2}-
\frac{R^2}{1+R^2}.
\end{align*}
We complete the proof of Proposition \ref{l2} by sending $\epsilon\to 0$
and $R\to\infty$.

Now we are ready to prove Theorem \ref{t1}. Suppose the infimum of
$I(\cdot,\Omega)$ is not attained, then Proposition \ref{l2} says
$E(\Omega)\ge -1-4\pi\gamma(\Omega)$. Together with Proposition \ref{l1}
we have
$$
E(\Omega)=-1-4\pi\gamma(\Omega).
$$
Since the infimum of $I(\cdot,B_1)$ is not attained,
$E(B_1)= -1-4\pi\gamma(B_1)=-1$. Therefore Lemma
\ref{gamma0} implies $E(\Omega)\ge E(B_1)$ and equality holds if
and only if $\Omega=B_1$.

If on the other hand, the infimum of $I(\cdot, \Omega)$ is attained by some
$u\in H^1_0(\Omega)$. Then we have $\Omega\ne B_1$ and
$$
E(\Omega)=I(u,\Omega)\ge I(u^*,B_1)> E(B_1),
$$
where $u^*:B_1\to \mathbf R$ is the rearrangement of $u$. This completes
the proof of Theorem \ref{t1}.
\medskip
\section{An existence result}
In this section we provide an existence result.
Let $\Omega\subset{\mathbf R}^2$ be a bounded domain with
$|\Omega|=\pi$. Suppose $\Omega$ can be covered by a strip $D_d$
with width $d\le \frac{\pi}{2\sqrt e}$, then we will show that the
infimum of $I(\cdot,\Omega)$ can be achieved.

Without loss of generality, we may assume that $D_d=\{z\ |\ 0<\mbox{Im}z<d\}$.
Here we have identified $\mathbf C$ with ${\mathbf R}^2$. It is easy to see
that $w=\phi(z)=e^{\pi z/d}$ maps $D_d$ to the upper half plane $\mbox{Im}z>0$.
Therefore for $\alpha\in(0,d)$, the Green's function of $D_d$ with
pole at $\alpha i$ is
$$
G(z,\alpha i)=\frac 1{2\pi}\log\left|
\frac{e^{\pi z/d}-e^{-\alpha \pi i/d}}{e^{\pi z/d}-e^{\alpha \pi i/d}}
\right|,
$$
and the regular part of the Green's function is given by
$$
\gamma(z,\alpha i)=\frac{1}{2\pi}\log\left|e^{\pi z/d}
-e^{-\alpha\pi i/d}\right|+\frac1{2\pi}\log\left|
\frac{z-\alpha i}{e^{\pi z/d}-e^{\alpha \pi i/d}}
\right|.
$$
Letting $z\to \alpha i$, we can easily see that
$$
\gamma(\alpha i)=\lim_{z\to\alpha i}\gamma(z,\alpha i)=\frac1{2\pi}
\log\frac{2\sin(\alpha \pi/d)}{\pi/d}\le\frac{1}{2\pi}\log\frac{2d}{\pi}.
$$
\begin{lemma}
Suppose $z_0\in\Omega_1\subset\Omega_2$, then $\gamma_{\Omega_1}(z_0)\le
\gamma_{\Omega_2}(z_0)$.
\end{lemma}

\begin{proof}
Let $G_{\Omega_1}(z,z_0)$ and $G_{\Omega_2}(z,z_0)$ be the Green's functions
on $\Omega_1$ and $\Omega_2$ respectively. Then the lemma follows from
applying maximum principle to $G_{\Omega_2}-G_{\Omega_1}$ on $\Omega_1$.
\end{proof}
Since $\Omega\subset D_d$, we have
\begin{equation}
\label{s1}
\gamma(\Omega)=\sup_{z}\gamma(z,\Omega)\le
\gamma(D_d)=\frac1{2\pi}\log\frac{2d}{\pi}.
\end{equation}
Suppose the infimum of $I(\cdot,\Omega)$ can not be achieved.
From Proposition \ref{l1} and Proposition \ref{l2} we have
$$
E(\Omega)=-1-4\pi \gamma(\Omega).
$$
On the other hand $E(\Omega)\le I(0,\Omega)=0$. Hence $\gamma(\Omega)
\ge -\frac1{4\pi}$. It follows from (\ref{s1}) that
$$
\frac{1}{2\pi}\log(2d/\pi)\ge -\frac 1{4\pi},
$$
or $d\ge\pi/(2\sqrt e)$, which contradicts with the assumption that
$d<\pi/(2\sqrt e)$. Hence the infimum of $I(\cdot,\Omega)$ is
attained.\\
{\bf Acknowledgement}:
I am very grateful to M. Zhu for valuable suggestions and
useful conversations.

\end{document}